# CONFIDENCE REGIONS FOR HIGH QUANTILES OF A HEAVY TAILED DISTRIBUTION

By Liang Peng[1] and Yongcheng Qi

*Georgia Institute of Technology and University of Minnesota Duluth*

Estimating high quantiles plays an important role in the context of risk management. This involves extrapolation of an unknown distribution function. In this paper we propose three methods, namely, the normal approximation method, the likelihood ratio method and the data tilting method, to construct confidence regions for high quantiles of a heavy tailed distribution. A simulation study prefers the data tilting method.

**1. Introduction.** In estimating high quantiles of an unknown probability distribution function, one has to infer beyond the observations. This can be done via extrapolating from intermediate quantiles when the underlying distribution has a regularly varying tail. An important application of high quantiles is to forecast rare events. Some references on this topic include [1, 3, 7, 14, 23, 25]. Like tail index estimation, only a part of the upper order statistics is involved in the estimation of high quantiles. Recently, Ferreira, de Haan and Peng [8] provided a data-driven method to choose the optimal sample fraction in terms of asymptotic mean squared errors.

In this paper we are interested in obtaining confidence regions for high quantiles. More specifically, three methods will be investigated, namely, the normal approximation method, the likelihood ratio method and the data tilting method; see Section 2 for details. We demonstrate by a simulation study that the data tilting method is preferred. Our data tilting method is similar to the general data tilting method proposed by Hall and Yao [12], which is employed to tilt time series data. This general data tilting method was applied to interval estimation, robust inference and inference under constraints for linear time series. One of its advantages is that it admits a

Received March 2004; revised September 2005.
[1] Supported in part by NSF Grant DMS-04-03443 and a Humboldt research fellowship.

wide range of distance functions. Tilting methods to statistics have a long history; nonparametric techniques involving tilting go back at least to work of Grenander [10], which studies nonparametric density estimation under monotonicity constraints.

To our best knowledge, not much work has been done in applying data tilting methods or empirical likelihood methods to statistics of extremes. The empirical likelihood method, introduced in [17, 18], is a nonparametric approach for constructing confidence regions. Like the bootstrap and the jackknife, the empirical likelihood method does not need to specify a family of distributions for the data. One of the advantages of the empirical likelihood method is that it enables the shape of a region, such as the degree of asymmetry in a confidence interval, to be determined automatically by the sample. In certain regular cases, empirical likelihood based confidence regions are Bartlett correctable; see [6, 11]. For a more complete disclosure of recent references and development, we refer to the book by Owen [19]. Recently, Lu and Peng [15] applied the empirical likelihood method to obtain confidence intervals for the tail index, and Peng [20] generalized the empirical likelihood method to the case of infinite variance. Here, we propose to employ the general data tilting method of Hall and Yao [12] to obtain confidence regions for high quantiles.

We organize this paper as follows. In Section 2 three different methods for constructing confidence regions for high quantiles are introduced, and main results about asymptotic limits are also presented. In Section 3 simulation results are reported for comparisons of the performance of the three methods in terms of both coverage probability and approximate interval length, and a real data application is also included. Finally, all the proofs are given in the Appendix.

**2. Methodologies and main results.** Let $X_1, \ldots, X_n$ be independent random variables with a common distribution function $F$ which satisfies

$$(1) \qquad 1 - F(x) = e(x)x^{-\gamma} \qquad \text{for } x > 0,$$

where $\gamma > 0$ is an unknown parameter called the tail index, and $e(x)$ is a slowly varying function, that is, $\lim_{t \to \infty} e(tx)/e(t) = 1$ for all $x > 0$. Let $X_{n,1} \leq \cdots \leq X_{n,n}$ denote the order statistics of $X_1, \ldots, X_n$.

Throughout this paper we assume that $p_n \in (0, 1)$ and $p_n \to 0$ as $n \to \infty$. A $100(1 - p_n)\%$ quantile for the distribution $F$ is defined as $x_p = (1 - F)^-(p_n)$, where $(\cdot)^-$ denotes the inverse function of $(\cdot)$. The main aim of this paper is to obtain confidence regions for $x_p$.

In order to introduce our methodologies, let us assume temporarily that $F$ has the simpler form

$$(2) \qquad 1 - F(x) = cx^{-\gamma} \qquad \text{for } x > T.$$

Put $\delta_i = I(X_i > T)$. Then the likelihood function for the censored data $\{(\delta_i, \max(X_i, T))\}_{i=1}^n$ is

$$L(\gamma, c) = \prod_{i=1}^n (c\gamma X_i^{-\gamma-1})^{\delta_i}(1 - cT^{-\gamma})^{1-\delta_i}.$$

In the paper we actually take $T = X_{n,n-k}$, where $k = k(n)$ satisfies

(3) $$k \to \infty \quad \text{and} \quad \frac{k}{n} \to 0.$$

Then the likelihood function above becomes

(4) $$L(\gamma, c) = \prod_{i=1}^n (c\gamma X_i^{-\gamma-1})^{\delta_i}(1 - cX_{n,n-k}^{-\gamma})^{1-\delta_i}.$$

Next we are ready to present our three methods for constructing confidence regions for $x_p$.

*Method* I: *Normal approximation method.* Let $(\hat{\gamma}_n, \hat{c}_n)$ denote the maximum likelihood estimator of $(c, \gamma)$, that is, $L(\hat{\gamma}_n, \hat{c}_n) = \max_{\gamma > 0, c > 0} L(\gamma, c)$. Then it is easy to check that $\hat{c}_n = \frac{k}{n} X_{n,n-k}^{\hat{\gamma}_n}$ and

$$\hat{\gamma}_n = \left\{ \frac{1}{k} \sum_{i=1}^k \log(X_{n,n-i+1}/X_{n,n-k}) \right\}^{-1}.$$

Note that $\hat{\gamma}_n$ is the well-known Hill estimator [13]. Therefore, by (2), a natural estimator for $x_p$ is $\hat{x}_p = (p_n/\hat{c}_n)^{-1/\hat{\gamma}_n}$. In order to derive the asymptotic normality of $\hat{x}_p$, we need a stricter condition than (1): suppose there exists a function $A(t) \to 0$ (as $t \to \infty$) such that

(5) $$\lim_{t \to \infty} \frac{U(tx)/U(t) - x^{1/\gamma}}{A(t)} = x^{1/\gamma} \frac{x^\rho - 1}{\rho}$$

for some $\rho < 0$, where $U(x) = (\frac{1}{1-F})^{\leftarrow}(x)$. Then $A(t)$ is a regularly varying function with index $\rho$; see [5]. Note that (5) implies (1). The following theorem can be derived from [8].

THEOREM 1. *Assume* (5) *and* (3) *hold. If* $\sqrt{k}A(n/k) \to 0$, $np_n = O(k)$ *and* $\log(\frac{k}{np_n})/\sqrt{k} \to 0$ *as* $n \to \infty$, *then*

(6) $$\frac{\hat{\gamma}_n \sqrt{k}}{\log(k/(np_n))} \log \frac{\hat{x}_p}{x_p} \xrightarrow{d} N(0, 1).$$

Hence, based on the above limit, a confidence interval with level $\alpha$ for $x_p$ is

$$I_\alpha^n = \left(\hat{x}_p \exp\left\{-z_\alpha \log\left(\frac{k}{np_n}\right) \Big/ (\hat{\gamma}_n \sqrt{k})\right\}, \hat{x}_p \exp\left\{z_\alpha \log\left(\frac{k}{np_n}\right) \Big/ (\hat{\gamma}_n \sqrt{k})\right\}\right),$$

where $z_\alpha$ satisfies $P(|N(0,1)| \leq z_\alpha) = \alpha$. This confidence interval has asymptotically correct coverage probability $\alpha$, that is, $P(x_p \in I_\alpha^n) \to \alpha$ as $n \to \infty$.

The next theorem presents the coverage expansion for $I_\alpha^n$.

THEOREM 2. *Under the conditions of Theorem 1,*

$$P\left(\frac{\hat{\gamma}_n \sqrt{k}}{\log(k/(np_n))} \log \frac{\hat{x}_p}{x_p} \leq x\right) - \Phi(x)$$

$$= \frac{1}{3\sqrt{k}} \phi(x)(1 + 2x^2) - \phi(x) \frac{\gamma}{1-\rho} \sqrt{k} A(n/k) - \frac{1}{2} x \phi(x) \left(\log \frac{k}{np_n}\right)^{-2}$$

$$+ o\left(\left(\log \frac{k}{np_n}\right)^{-2} + \frac{1}{\sqrt{k}} + \sqrt{k} |A(n/k)|\right),$$

*uniformly for* $-\infty < x < \infty$, *where* $\Phi(x)$ *and* $\phi(x)$ *denote the distribution function and density function of* $N(0,1)$, *respectively. Furthermore,*

$$P(x_p \in I_\alpha^n) = \alpha - z_\alpha \phi(z_\alpha) \left(\log \frac{k}{np_n}\right)^{-2}$$

$$+ o\left(\left(\log \frac{k}{np_n}\right)^{-2} + \frac{1}{\sqrt{k}} + \sqrt{k} |A(n/k)|\right).$$

REMARK 1. Theorem 2 shows that $P(x_p \in I_\alpha^n) - \alpha = O((\log n)^{-2})$ in the case $\log(np_n) = O(\log(n))$. This means that the coverage accuracy for high quantiles is not very accurate in general. To achieve this asymptotic rate, $k$ can be of order $n^\theta$ for some $0 < \theta < -2\rho/(1-2\rho)$. The unknown parameter $\rho$ can be estimated; see, for example, [21].

*Method* II: *Likelihood ratio method.* Define $\hat{\gamma}_n$ and $\hat{c}_n$ as in Method I. First set

$$l_1 = \max_{\gamma > 0, c > 0} \log L(\gamma, c) = \log L(\hat{\gamma}_n, \hat{c}_n).$$

Next we maximize $\log L(\gamma, c)$ subject to

$$\gamma > 0, \qquad c > 0, \qquad \gamma \log x_p + \log\left(\frac{p_n}{c}\right) = 0,$$

and denote this maximized likelihood function by $l_2(x_p)$. Note that the above equation comes from setting $p_n = 1 - F(x_p) = c x_p^{-\gamma}$. It is easy to show that

$$l_2(x_p) = \log L(\bar{\gamma}(\lambda), \bar{c}(\lambda)),$$

where

$$\bar{\gamma}(\lambda) = \frac{k}{\sum_{i=1}^{k}(\log X_{n,n-i+1} - \log X_{n,n-k}) + \lambda \log X_{n,n-k} - \lambda \log x_p},$$

$$\bar{c}(\lambda) = X_{n,n-k}^{\bar{\gamma}(\lambda)} \frac{k-\lambda}{n-\lambda}$$

and $\lambda$ satisfies

(7) $$\bar{\gamma}(\lambda) \log x_p + \log\left(\frac{p_n}{\bar{c}(\lambda)}\right) = 0,$$

(8) $$\bar{\gamma}(\lambda) > 0 \quad \text{and} \quad \lambda < k.$$

Therefore, the log-likelihood ratio multiplied by minus two is

$$l(x_p) = -2(l_2(x_p) - l_1).$$

THEOREM 3. *Suppose the conditions in Theorem 1 hold. Then there exists a unique solution to (7) and (8), say, $\hat{\lambda}(x_p)$, and*

(9) $$l(x_{p,0}) \xrightarrow{d} \chi^2(1),$$

*with $\lambda = \hat{\lambda}(x_{p,0})$ in the definition of $l_2(x_{p,0})$, where $x_{p,0}$ is the true value of $x_p$.*

Therefore, based on the above limit, a confidence region with level $\alpha$ for $x_p$ is

$$I_\alpha^l = \{x_p : l(x_p) \leq u_\alpha\},$$

where $u_\alpha$ is the $\alpha$-level critical point of $\chi^2(1)$. This confidence region has asymptotically correct coverage probability $\alpha$, that is, $P(x_p \in I_\alpha^l) \to \alpha$ as $n \to \infty$.

REMARK 2. The profile likelihood approach has been employed to construct confidence regions for high quantiles based on fitting a generalized Pareto distribution to exceedances over a deterministic high threshold; see [24]. The difference between our Method II and the profile likelihood method is that we take the random high threshold into account in our censored likelihood function.

*Method* III: *Data tilting method.* Here we employ a data tilting method, similar to that of Hall and Yao [12], to construct a confidence region for $x_p$. First, for any fixed weights $q = (q_1, \ldots, q_n)$ such that $q_i \geq 0$ and $\sum_{i=1}^n q_i = 1$, we solve

$$(\hat{\gamma}(q), \hat{c}(q)) = \arg\max_{(\gamma,c)} \sum_{i=1}^n q_i \log((c\gamma X_i^{-\gamma-1})^{\delta_i}(1 - cX_{n,n-k}^{-\gamma})^{1-\delta_i}).$$

This results in
$$\hat{\gamma}(q) = \frac{\sum_{i=1}^n q_i \delta_i}{\sum_{i=1}^n q_i \delta_i (\log X_i - \log X_{n,n-k})},$$
$$\hat{c}(q) = X_{n,n-k}^{\hat{\gamma}(q)} \sum_{i=1}^n q_i \delta_i.$$

Define
$$D_{\rho_0}(q) = \begin{cases} (\rho_0(1-\rho_0))^{-1}\left(1 - n^{-1}\sum_{i=1}^n (nq_i)^{\rho_0}\right), & \text{if } \rho_0 \neq 0, 1, \\ -n^{-1}\sum_{i=1}^n \log(nq_i), & \text{if } \rho_0 = 0, \\ \sum_{i=1}^n q_i \log(nq_i), & \text{if } \rho_0 = 1. \end{cases}$$

The function $D_{\rho_0}(q)$ is a measure of distance between $q$ and uniform distribution, that is, $q_i = 1/n$. Next, we shall choose $q$ to minimize this distance. More specifically, solve $(2n)^{-1}L(x_p) = \min_q D_{\rho_0}(q)$ subject to the constraints
$$q_i \geq 0, \qquad \sum_{i=1}^n q_i = 1, \qquad \hat{\gamma}(q)\log \frac{x_p}{X_{n,n-k}} = \log \frac{\sum_{i=1}^n q_i \delta_i}{p_n}.$$

The constraint $\hat{\gamma}(q)\log(x_p/X_{n,n-k}) = \log(\sum_{i=1}^n q_i \delta_i/p_n)$ is equivalent to $x_p = (p_n/\hat{c}(q))^{-1/\hat{\gamma}(q)}$.

Here we only consider the case $\rho_0 = 1$ since other cases are similar and the case $\rho_0 = 1$ gives good robustness properties. Put
$$A_1(\lambda_1) = 1 - \frac{n-k}{n}e^{-1-\lambda_1} \quad \text{and} \quad A_2(\lambda_1) = A_1(\lambda_1)\frac{\log(x_p/X_{n,n-k})}{\log(A_1(\lambda_1)/p_n)}.$$

Then, by the standard method of Lagrange multipliers, we have

(10)
$$q_i = q_i(\lambda_1, \lambda_2) = \begin{cases} \frac{1}{n}e^{-1-\lambda_1}, & \text{if } \delta_i = 0, \\ \frac{1}{n}\exp\Big\{-1-\lambda_1 \\ \qquad + \lambda_2\Big(\frac{\log(x_p/X_{n,n-k})}{A_2(\lambda_1)} - \frac{1}{A_1(\lambda_1)} \\ \qquad - \frac{A_1(\lambda_1)\log(X_i/X_{n,n-k})\log(x_p/X_{n,n-k})}{A_2^2(\lambda_1)}\Big)\Big\}, \\ & \text{if } \delta_i = 1, \end{cases}$$

where $\lambda_1$ and $\lambda_2$ satisfy

(11) $$\sum_{i=1}^n q_i = 1, \qquad \hat{\gamma}(q)\log \frac{x_p}{X_{n,n-k}} = \log \frac{\sum_{i=1}^n q_i \delta_i}{p_n}.$$

THEOREM 4. *Suppose the conditions in Theorem 1 hold. Then, with probability tending to one, there exists a solution to (11), say, $(\hat{\lambda}_1(x_p), \hat{\lambda}_2(x_p))$, such that, for $(\lambda_1, \lambda_2) = (\hat{\lambda}_1(x_p), \hat{\lambda}_2(x_p))$,*

$$(12) \quad -\log\left(1 + \frac{\sqrt{k}\sqrt{\log(k/(np_n))}}{n-k}\right) \leq 1 + \lambda_1$$

$$\leq -\log\left(1 - \frac{\sqrt{k}\sqrt{\log(k/(np_n))}}{n-k}\right),$$

$$(13) \quad |\lambda_2| \leq k^{-1/4}\frac{k/n}{\log(k/(np_n))}$$

*and $L(x_{p,0}) \xrightarrow{d} \chi^2(1)$ with $(\lambda_1, \lambda_2) = (\hat{\lambda}_1(x_{p,0}), \hat{\lambda}_2(x_{p,0}))$ in the definition of $L(x_{p,0})$.*

Hence, based on the above limit, a confidence region with level $\alpha$ for $x_p$ is

$$I_\alpha^t = \{x_p : L(x_p) \leq u_\alpha\},$$

where $u_\alpha$ is the $\alpha$-level critical point of $\chi^2(1)$. This confidence region has asymptotically correct coverage probability $\alpha$, that is, $P(x_p \in I_\alpha^t) \to \alpha$ as $n \to \infty$.

REMARK 3. In order to compare these three methods theoretically, it is necessary to derive corresponding coverage expansions for $I_\alpha^l$ and $I_\alpha^t$. This requires much work and it will be one of our future topics.

### 3. Simulation study and real application.

3.1. *A simulation study.* In order to compare the performance of confidence regions based on the normal approximation method, the likelihood ratio method and the data tilting method, we conducted a simulation study to examine coverage probabilities and approximate lengths of confidence regions.

We employed the following two distributions: (i) the Burr$(\alpha, \beta)$ distribution, given by $F(x) = 1 - (1 + x^{\beta-\alpha})^{-\alpha/(\beta-\alpha)}$ $(x > 0)$; (ii) the Fréchet$(\alpha)$ distribution, given by $F(x) = \exp(-x^{-\alpha})$ $(x > 0)$. Corresponding to (5), we have $\gamma = 1/\alpha$, $\rho = -\frac{\beta-\alpha}{\alpha}$ and $A(t) = \alpha^{-1}t^{-(\beta-\alpha)/\alpha}$ for Burr$(\alpha, \beta)$, and $\gamma = 1/\alpha$, $\rho = -1$ and $A(t) = (2\alpha t)^{-1}$ for Fréchet$(\alpha)$.

First, we generated 10,000 random samples of size $n = 1000$ from the distributions Burr$(1, 1.5)$, Burr$(1, 2)$, Burr$(2, 3)$, Burr$(2, 4)$, Fréchet$(1)$ and Fréchet$(2)$, and then computed coverage probabilities of $I_{0.9}^n, I_{0.9}^l$ and $I_{0.9}^t$ for $p_n = 0.01$ and $p_n = 0.001$. These coverage probabilities are plotted against different sample fractions $k = 20, 25, \ldots, 300$ in Figures 1–4. From these figures we observe that the data tilting method is better than the other two methods in terms of coverage accuracy in general, especially for larger values of $k$. Thus, the data tilting method is less sensitive to the bias when a large value of $k$ is employed. This may be due to the automatic choice of weights $q_i$ in the data tilting method. Although it does not make much sense to compare these three methods with the empirical likelihood method for quantiles (see Section 3.6 of [19]), we find that the coverage probabilities based on the empirical likelihood method for quantiles are 0.7631 for $p_n = 0.01$ and 0.6047 for $p_n = 0.001$. These coverage probabilities are not as accurate as those based on the other three methods for most sample fractions $k$. Note that the empirical likelihood method for quantiles is independent of the underlying distribution function.

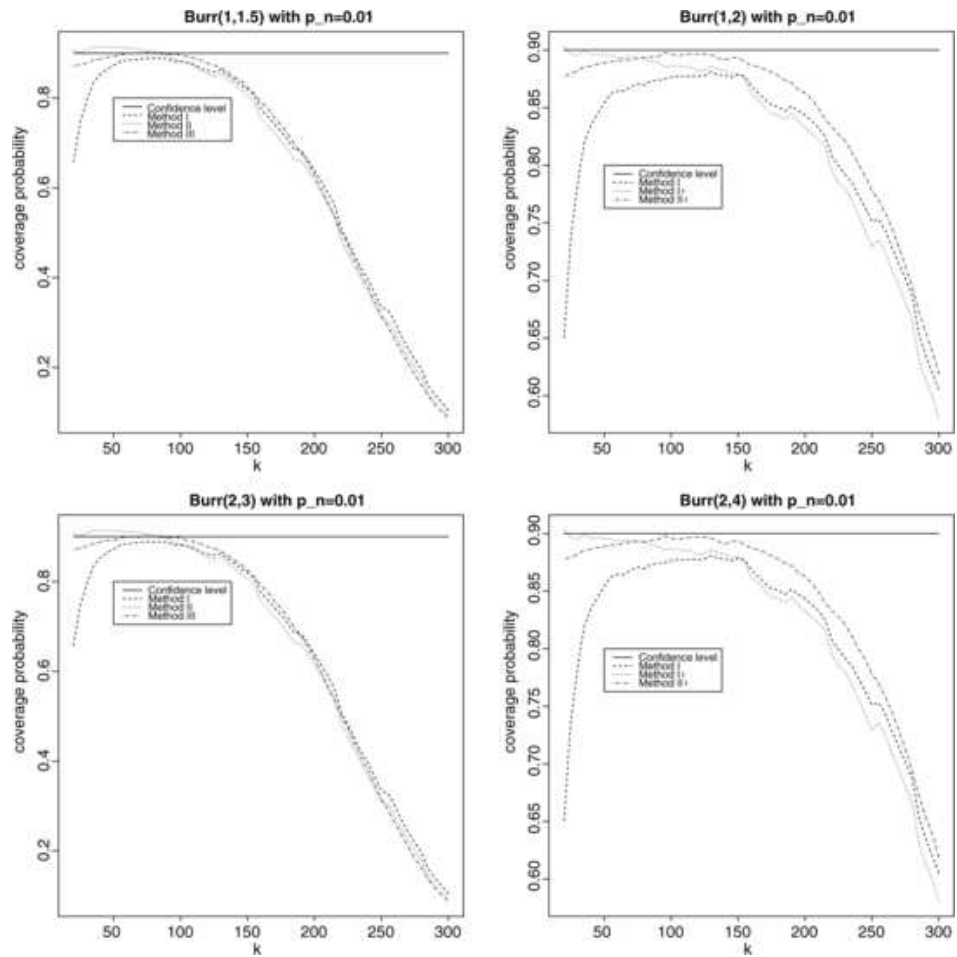

Fig. 1. *Coverage probabilities for Burr distributions with $p_n = 0.01$. The coverage probabilities of confidence regions $I_{0.90}^n$, $I_{0.90}^l$ and $I_{0.90}^t$ are plotted against the different sample fractions $k = 20, 25, \ldots, 300$ for different Burr distributions.*

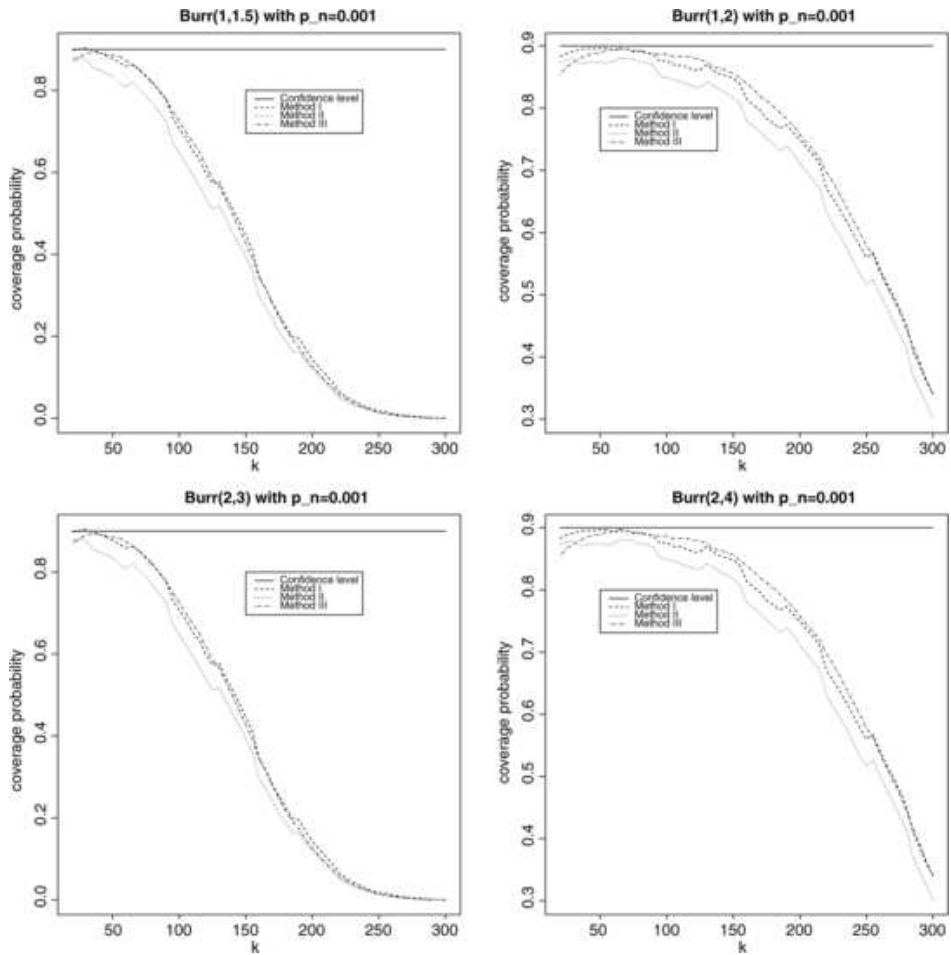

Fig. 2. *Coverage probabilities for Burr distributions with $p_n = 0.001$. The coverage probabilities of confidence regions $I^n_{0.90}$, $I^l_{0.90}$ and $I^t_{0.90}$ are plotted against the different sample fractions $k = 20, 25, \ldots, 300$ for different Burr distributions.*

Second, we generated 1,000 random samples of size $n = 1000$ from the distributions Burr$(1, 2)$ and Fréchet$(1)$, and then calculated the length of $I^n_{0.9}$ and the approximate lengths of $I^l_{0.9}$ and $I^t_{0.9}$ for $p_n = 0.01$. Let us explain how we calculate the approximate length of $I^t_{0.9}$. The same algorithm was employed to obtain the approximate length of $I^l_{0.9}$. First, we search an $x_p$ near $\hat{x}_p$ such that $L(x_p) < u_{0.9}$. Then we both increase and decrease $x_p$ by a small step 0.1 until $L(x_p) > u_{0.9}$. The corresponding values are denoted by $x^u_p$ and $x^l_p$, respectively. Thus, we approximate $I^t_{0.9}$ by the interval $[x^l_p, x^u_p]$.

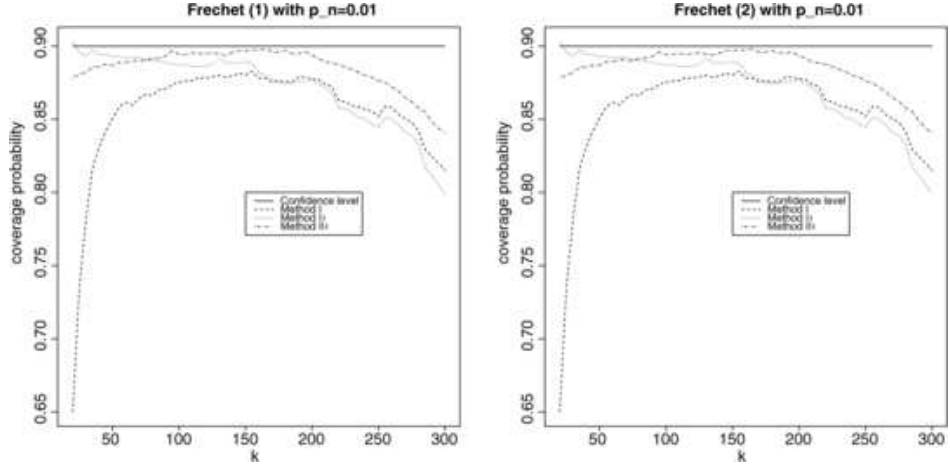

Fig. 3. *Coverage probabilities for Fréchet distributions with $p_n = 0.01$. The coverage probabilities of confidence regions $I_{0.90}^n$, $I_{0.90}^l$ and $I_{0.90}^t$ are plotted against the different sample fractions $k = 20, 25, \ldots, 300$ for different Fréchet distributions.*

These approximate lengths, $x_p^u - x_p^l$, are plotted against the different sample fractions $k = 20, 30, \ldots, 300$ in Figure 5. We notice that the approximate confidence interval lengths based on the data tilting method are smallest for most cases.

Third, we generated a random sample of size $n = 1000$ from the distributions Burr$(1, 2)$ and Fréchet$(1)$, and then computed the data tilting like-lihood function $L(x_p)$ for $p_n = 0.01$ and $x_p = x_{p,0} - 50 + i$, $i = 0, 1, \ldots, 200$, where $x_{p,0}$ denotes the true quantile. We took $k = 50$ and $100$. Figure 6 indicates that the data tilting likelihood function is approximately convex, which suggests that $I_{0.9}^t$ may indeed be an interval. Unlike the empirical likelihood method for means, we were unable to prove that $I_\alpha^t$ is an interval.

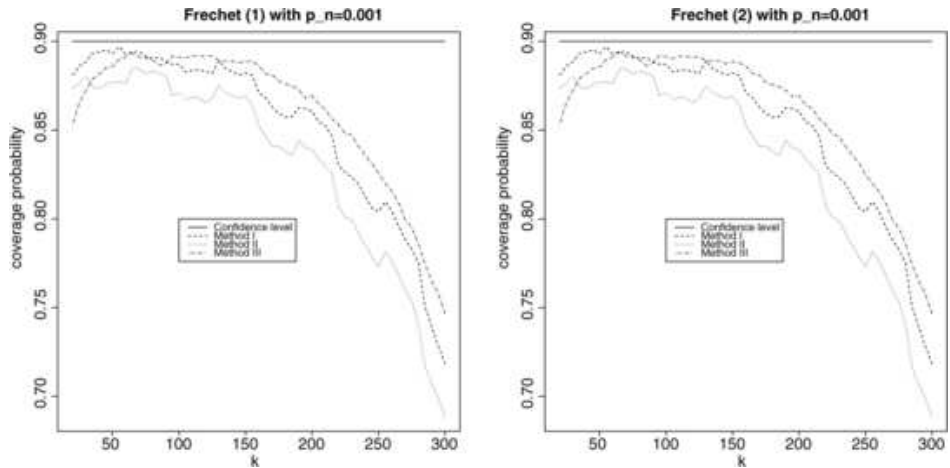

Fig. 4. *Coverage probabilities for Fréchet distributions with $p_n = 0.001$. The coverage probabilities of confidence regions $I_{0.90}^n$, $I_{0.90}^l$ and $I_{0.90}^t$ are plotted against the different sample fractions $k = 20, 25, \ldots, 300$ for different Fréchet distributions.*

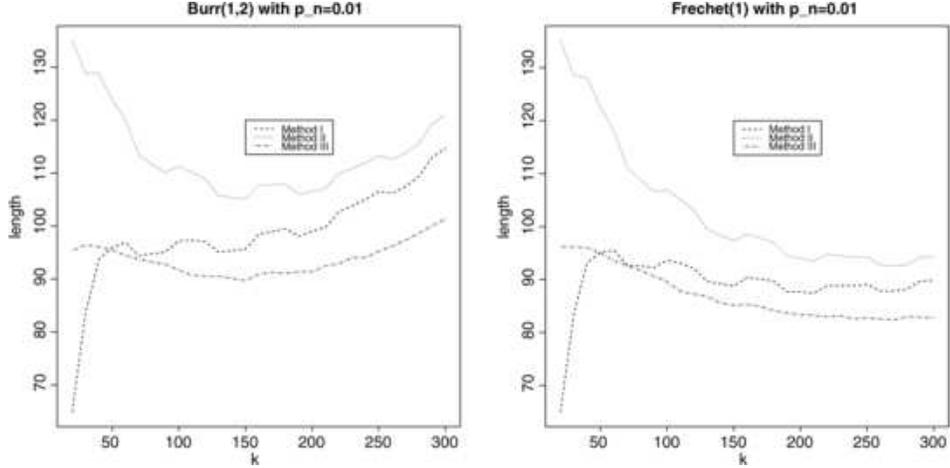

FIG. 5. *Averages of the approximate confidence lengths with $p_n = 0.01$. The averages of approximate lengths of confidence regions $I_{0.90}^n$, $I_{0.90}^l$ and $I_{0.90}^t$ are plotted against the different sample fractions $k = 20, 30, \ldots, 300$ for* Burr$(1.0, 2.0)$ *and Fréchet$(1)$ distributions.*

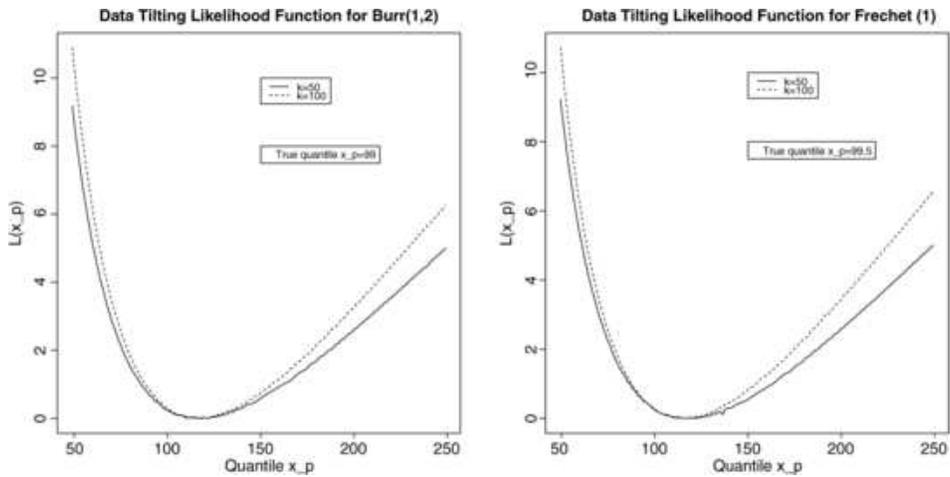

FIG. 6. *Data tilting likelihood function $L(x_p)$ with $p_n = 0.01$. The data tilting likelihood functions are plotted against different $x_p = x_{p,0} - 50 + i$, $i = 0, 1, \ldots, 200$, for* Burr$(1.0, 2.0)$ *and Fréchet$(1)$, where $x_{p,0}$ denotes the true quantile. We took $k = 50$ and $k = 100$.*

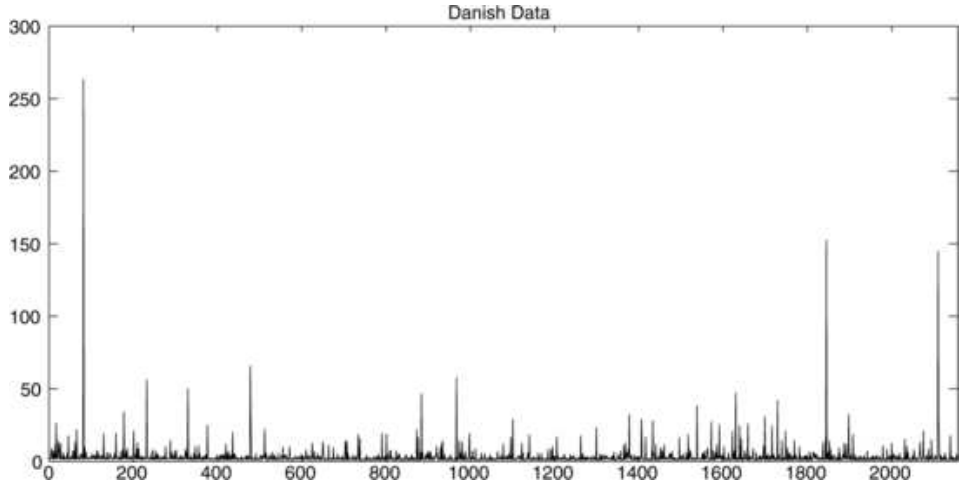

Fig. 7. *Danish fire loss data. This consists of* 2156 *losses over one million Danish Kroner (DKK) from the years* 1980 *to* 1990, *inclusive.*

In summary, our simulation study for sample size $n = 1000$ prefers the data tilting method, which gives the best coverage accuracy in general, is less sensitive to the choice of sample fraction $k$, and has a shorter approximate interval length in most cases. Although we do not report the simulation study for sample size $n = 200$, the same conclusions as above are drawn, except that Method II performs worst.

3.2. *A real application.* The data set we shall analyze consists of 2156 Danish fire losses over one million Danish Kroner (DKK) from the years 1980 to 1990 inclusive (see Figure 7). The loss figure is a total loss for the event concerned and includes damage to buildings, damage to furnishings and personal property, as well as loss of profits. This data set was analyzed by McNeil [16] and Resnick [22], where the right tail index was confirmed to be between 1 and 2. Further, Peng [20] applied the empirical likelihood method to this data set to obtain a confidence interval for the mean.

We took $p_n = 0.001$ and plotted the confidence interval $I_{0.90}^n$, and the approximate confidence intervals $I_{0.90}^l$ and $I_{0.90}^t$, against the different sample fraction $k = 60, 65, \ldots, 400$ in Figure 8. We note again that the approximate interval lengths based on the data tilting method are smallest for most cases.

## APPENDIX A: PROOFS OF THEOREMS 2, 3 AND 4

PROOF OF THEOREM 2. Let $V_1, V_2, \ldots, V_n$ be i.i.d. random variables uniformly distributed over $(0, 1)$ and $V_{n,1} \leq V_{n,2} \leq \cdots \leq V_{n,n}$ be the order statistics of $V_1, V_2, \ldots, V_n$. Define $c_n = 1 - \frac{k}{n+1}$ and $d_n = \sqrt{c_n(1-c_n)/(n+1)}$.

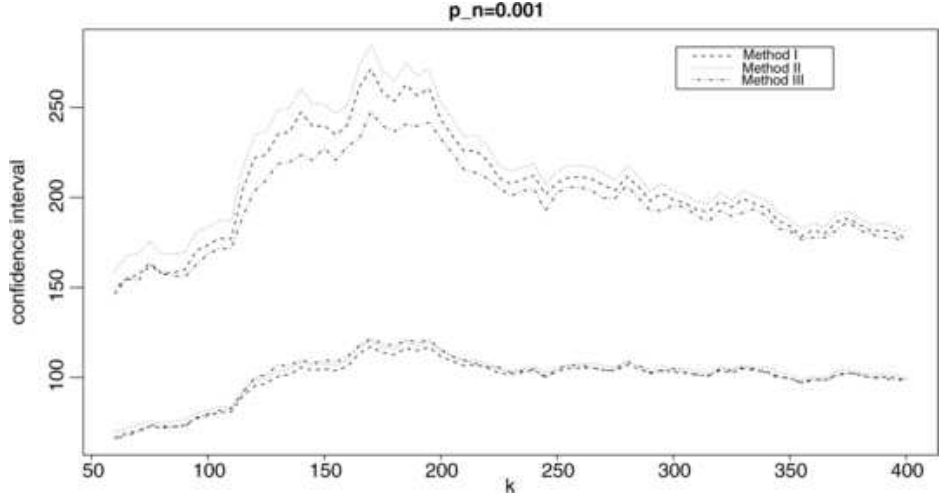

Fig. 8. *Approximate confidence intervals for Danish fire loss data. The approximate confidence intervals with level* 0.90 *based on the normal approximation method (Method* I*), the likelihood ratio method (Method* II*) and the data tilting method (Method* III*) are plotted against* $k = 60, 65, \ldots, 400$. *We took* $p_n = 0.001$.

Then the density function of $(V_{n,n-k} - c_n)/d_n$ is

$$\phi_n(u) = \begin{cases} \dfrac{n! d_n}{k!(n-k-1)!}(c_n + d_n u)^{n-k-1}(1 - c_n - d_n u)^k, \\ \qquad \text{if } 0 < c_n + d_n u < 1, \\ 0, \qquad \text{otherwise.} \end{cases}$$

For each $-c_n/d_n < u < (1-c_n)/d_n$, define

$$Y_{n,j}(u) = \gamma\{\log U((1-c_n-d_n u)^{-1}(1-V_i)^{-1}) - \log U((1-c_n-d_n u)^{-1})\},$$

$$H_n = \sqrt{k}\gamma(\hat\gamma_n^{-1} - \gamma^{-1}),$$

$$H_n(u) = \frac{1}{\sqrt{k}}\sum_{j=1}^k (Y_{n,j}(u) - 1),$$

$$r_n(u) = \sqrt{k}\left\{\frac{\gamma}{\log(k/(np_n))}\log\frac{U(p_n^{-1})}{U(1/(1-c_n-d_n u))} - 1\right\}.$$

Since

$$\frac{\hat\gamma_n\sqrt{k}}{\log(k/(np_n))}\log\frac{\hat x_p}{x_p} = \frac{\hat\gamma_n\sqrt{k}}{\log(k/(np_n))}\log\frac{X_{n,n-k}}{x_p} + \sqrt{k}$$

and

$$P\left(\frac{\hat\gamma_n\sqrt{k}}{\log(k/(np_n))}\log\frac{\hat x_p}{x_p} \leq x\right)$$

$$= P\left(\hat{\gamma}_n^{-1} \leq \frac{1}{1 - x/\sqrt{k}} \frac{1}{\log(k/(np_n))} \log \frac{x_p}{X_{n,n-k}}\right)$$

$$= P\left(H_n \leq \frac{1}{1 - x/\sqrt{k}} \left\{x + \sqrt{k}\left(\frac{\gamma}{\log(k/(np_n))} \log \frac{x_p}{X_{n,n-k}} - 1\right)\right\}\right)$$

for $|x| \leq k^{1/4}$, it follows from Lemma 2.2 of [2] that

$$(14) \quad P\left(\frac{\hat{\gamma}_n \sqrt{k}}{\log(k/(np_n))} \log \frac{\hat{x}_p}{x_p} \leq x\right) = \int_{-\infty}^{\infty} P\left(H_n(u) \leq \frac{x + r_n(u)}{1 - x/\sqrt{k}}\right) \phi_n(u) \, du$$

for $|x| \leq k^{-1/4}$. Similar to Lemma 2.3 of [2], we can prove that

$$P\left(H_n(u) \leq \frac{x + r_n(u)}{1 - x/\sqrt{k}}\right)$$

$$(15) \quad = \Phi\left(\frac{x + r_n(u)}{1 - x/\sqrt{k}}\right) + \phi\left(\frac{x + r_n(u)}{1 - x/\sqrt{k}}\right) \frac{1}{3\sqrt{k}} \left\{1 - \left(\frac{x + r_n(u)}{1 - x/\sqrt{k}}\right)^2\right\}$$

$$- \phi\left(\frac{x + r_n(u)}{1 - x/\sqrt{k}}\right) \frac{\gamma}{1 - \rho} \sqrt{k} A(n/k) + o\left(\frac{1}{\sqrt{k}} + \sqrt{k}|A(n/k)|\right),$$

uniformly for $x \in R$ and $|u| \leq k^{1/4}$. Since (5) is equivalent to

$$\lim_{t \to \infty} \frac{\log(U(tx)) - \log(U(t)) - \log(x)/\gamma}{A(t)} = \frac{x^\rho - 1}{\rho},$$

it follows from Potter's bounds that

$$r_n(u) = \frac{u}{\log(k/(np_n))} + \frac{u^2}{2\sqrt{k}\log(k/(np_n))} + \frac{u^3}{3k\log(k/(np_n))}$$

$$+ O\left(\frac{\sqrt{k}|A(n/k)|}{\log(k/(np_n))} + \frac{1}{\sqrt{k}\log(k/(np_n))}\right),$$

uniformly for $|u| \leq k^{1/4}$. Hence,

$$\Phi\left(\frac{x + r_n(u)}{1 - x/\sqrt{k}}\right) - \Phi(x)$$

$$(16) \quad = \phi(x)\frac{u}{\log(k/(np_n))} + \phi(x)x^2/\sqrt{k} - \frac{1}{2}x\phi(x)u^2\left(\log \frac{k}{np_n}\right)^{-2}$$

$$+ o\left(\left(\log \frac{k}{np_n}\right)^{-2} + \frac{1}{\sqrt{k}} + \frac{\sqrt{k}|A(n/k)|}{\log(k/(np_n))}\right)(1 + u^2),$$

uniformly for $|x| \leq k^{1/4}$ and $|u| \leq k^{1/4}$. Thus, the theorem follows from (14)–(16) and the facts that $\int |u|^t \phi_n(u) \, du = O(1)$ and $\int |u|^t \phi_n(u) I(|u| > k^{1/4}) \, du = o(1/k)$ for any $t > 0$. □

PROOF OF THEOREM 3. By the definition of $\bar{\gamma}(\lambda)$, we have

(17) $$\bar{\gamma}(\lambda) = \frac{\hat{\gamma}_n}{1 - (\lambda/k)\hat{\gamma}_n \log(x_p/X_{n,n-k})}.$$

Thus, equations (7) and (8) are equivalent to

(18) $$\frac{\hat{\gamma}_n \log(x_p/X_{n,n-k})}{1 - (\lambda/k)\hat{\gamma}_n \log(x_p/X_{n,n-k})} + \log \frac{(n-\lambda)p_n}{k-\lambda} = 0$$

and

(19) $$1 - \frac{\lambda}{k}\hat{\gamma}_n \log(x_p/X_{n,n-k}) > 0 \quad \text{and} \quad \lambda < k,$$

respectively. Set

$$g(\lambda) = \frac{\hat{\gamma}_n \log(x_p/X_{n,n-k})}{1 - (\lambda/k)\hat{\gamma}_n \log(x_p/X_{n,n-k})} + \log \frac{(n-\lambda)p_n}{k-\lambda}.$$

Then $g(\lambda)$ is continuous and increasing in $\lambda$ under the restriction (19) since

$$g'(\lambda) = \frac{[\hat{\gamma}_n \log(x_p/X_{n,n-k})]^2}{k[1 - (\lambda/k)\hat{\gamma}_n \log(x_p/X_{n,n-k})]^2} + \frac{n-k}{(n-\lambda)(k-\lambda)} > 0.$$

If $x_p/X_{n,n-k} > 1$, then (19) is equivalent to

$$\lambda < \min(k, k[\hat{\gamma}_n \log(x_p/X_{n,n-k})]^{-1}) =: a_n.$$

Since $g(-\infty) = \log p_n < 0$ and $g(a_n-) = \infty$, we conclude that there exists a unique $\lambda < a_n$ such that $g(\lambda) = 0$, that is, there exists a unique $\lambda$ satisfying (18) and (19). We can draw the same conclusion for the cases $x_p/X_{n,n-k} = 1$ and $x_p/X_{n,n-k} < 1$. So we prove the existence and uniqueness of a solution to (7) and (8).

Since

$$\log \hat{x}_p = -\frac{1}{\hat{\gamma}_n} \log \frac{p_n}{\hat{c}_n} = -\frac{1}{\hat{\gamma}_n}\left(\log p_n - \log \frac{k}{n} - \hat{\gamma}_n \log X_{n,n-k}\right),$$

we have

$$\log(\hat{x}_p/X_{n,n-k}) = -\frac{1}{\hat{\gamma}_n}\left(\log p_n - \log \frac{k}{n}\right) = \frac{1}{\hat{\gamma}_n} \log \frac{k}{np_n}$$

and

$$\hat{\gamma}_n \log(x_p/X_{n,n-k}) = \log \frac{k}{np_n} - \hat{\gamma}_n \log \frac{\hat{x}_p}{x_p}.$$

It follows from (6) that $\frac{\hat{x}_p}{x_p} - 1 \xrightarrow{p} 0$. Thus,

(20) $$\hat{\gamma}_n \log(x_p/X_{n,n-k}) = \left(\log \frac{k}{np_n}\right)(1 + o_p(1)) \xrightarrow{p} \infty.$$

Denote the solution to (7) and (8) by $\lambda_n$. Thus, $g(\lambda_n) = 0$. First we will show that

$$(21) \qquad P\left(|\lambda_n| < \frac{k}{[\hat{\gamma}_n \log(x_p/X_{n,n-k})]^2}\right) \to 1.$$

This is equivalent to proving that, as $n \to \infty$,

$$(22) \qquad P(g(b_n) > 0) \to 1 \quad \text{and} \quad P(g(-b_n) < 0) \to 1,$$

where $b_n = k/[\hat{\gamma}_n \log(x_p/X_{n,n-k})]^2$.

By (20) and Taylor's expansion,

$$g(b_n) = \hat{\gamma}_n \log(x_p/X_{n,n-k})\left(1 + \frac{b_n}{k}\hat{\gamma}_n \log(x_p/X_{n,n-k})(1 + o_p(1))\right)$$

$$- \log\frac{k}{np_n} + \log\left(1 - \frac{b_n}{n}\right) - \log\left(1 - \frac{b_n}{k}\right)$$

$$= \hat{\gamma}_n \log(x_p/X_{n,n-k}) - \log\frac{k}{np_n} + \frac{b_n}{k}(1 + o_p(1))$$

$$+ \frac{b_n}{k}[\hat{\gamma}_n \log(x_p/X_{n,n-k})]^2(1 + o_p(1))$$

$$= -\hat{\gamma}_n \log\frac{\hat{x}_p}{x_p} + \frac{b_n}{k}[\hat{\gamma}_n \log(x_p/X_{n,n-k})]^2(1 + o_p(1))$$

$$= 1 + o_p(1).$$

Similarly, $g(-b_n) = -1 + o_p(1)$. This yields (22) and, hence, (21).

Since (20) and (21) imply

$$(23) \qquad \frac{\lambda_n}{k}\hat{\gamma}_n \log(x_p/X_{n,n-k}) \xrightarrow{p} 0,$$

using Taylor's expansion again, we have

$$0 = g(\lambda_n) = -\hat{\gamma}_n \log\frac{\hat{x}_p}{x_p} + \frac{\lambda_n}{k}[\hat{\gamma}_n \log(x_p/X_{n,n-k})]^2(1 + o_p(1)),$$

that is,

$$(24) \quad \frac{\lambda_n}{k} = \frac{\hat{\gamma}_n \log(\hat{x}_p/x_p)}{[\hat{\gamma}_n \log(x_p/X_{n,n-k})]^2}(1 + o_p(1)) = \frac{\hat{\gamma}_n \log(\hat{x}_p/x_p)}{(\log(k/(np_n)))^2}(1 + o_p(1)).$$

Note that

$$\log L(\gamma, c) = k\log(c\gamma) - (\gamma + 1)\sum_{i=1}^{k} \log X_{n,n-i+1}$$

$$+ (n-k)\log(1 - cX_{n,n-k}^{-\gamma})$$

$$= k\log(\gamma) - (\gamma + 1)k\hat{\gamma}_n^{-1} + k\log(cX_{n,n-k}^{-\gamma}) - k\log(X_{n,n-k})$$

$$+ (n-k)\log(1 - cX_{n,n-k}^{-\gamma})$$

and

$$l(x_p) = -2(\log L(\bar{\gamma}(\lambda_n), \bar{c}(\lambda_n)) - \log L(\hat{\gamma}_n, \hat{c}_n))$$
$$= -2k\left(\log \frac{\bar{\gamma}(\lambda_n)}{\hat{\gamma}_n} - \left(\frac{\bar{\gamma}(\lambda_n)}{\hat{\gamma}_n} - 1\right)\right) - 2k\log\left(1 - \frac{\lambda_n}{k}\right) + 2n\log\left(1 - \frac{\lambda_n}{n}\right).$$

In view of (17) and (23), we have

$$\frac{\bar{\gamma}(\lambda_n)}{\hat{\gamma}_n} = \frac{1}{1 - (\lambda_n/k)\hat{\gamma}_n \log(x_p/X_{n,n-k})} \xrightarrow{p} 1.$$

It follows from (6) and (24) that

(25) $$\lambda_n/\sqrt{k} \xrightarrow{p} 0.$$

Hence, by (20), (6), (25) and Taylor's expansion,

$$\begin{aligned} l(x_p) &= k\left(\frac{\bar{\gamma}(\lambda_n)}{\hat{\gamma}_n} - 1\right)^2 (1 + o_p(1)) + O_p(\lambda_n^2/k) \\ &= \frac{(\lambda_n \hat{\gamma}_n \log(x_p/X_{n,n-k}))^2}{k}(1 + o_p(1)) + o_p(1) \\ &= \left(\frac{\hat{\gamma}_n \sqrt{k} \log(\hat{x}_p/x_p)}{\hat{\gamma}_n \log(x_p/X_{n,n-k})}\right)^2 (1 + o_p(1)) + o_p(1) \\ &= \left(\frac{\hat{\gamma}_n \sqrt{k} \log(\hat{x}_p)x_p}{\log(k/(np_n))}\right)^2 (1 + o_p(1)) + o_p(1) \\ &\xrightarrow{d} \chi^2(1). \quad \square \end{aligned}$$

PROOF OF THEOREM 4. Let $Z_i = \log(X_{n,n-i+1}/X_{n,n-k})$ for $1 \le i \le k$, $q_{(i)} = \frac{1}{n}\exp\{-1 - \lambda_1\}$ for $k + 1 \le i \le n$ and

$$q_{(i)} = \frac{1}{n}\exp\left\{-1 - \lambda_1 + \lambda_2\left(\frac{\log(x_p/X_{n,n-k})}{A_2(\lambda_1)} - \frac{1}{A_1(\lambda_1)}\right. \right.$$
$$\left.\left. - \frac{A_1(\lambda_1)Z_i \log(x_p/X_{n,n-k})}{A_2^2(\lambda_1)}\right)\right\}$$

for $1 \le i \le k$. Then (11) is equivalent to

$$\sum_{i=1}^k q_{(i)} = A_1(\lambda_1) \quad \text{and} \quad \sum_{i=1}^k q_{(i)} Z_i = A_2(\lambda_1).$$

Furthermore, this is equivalent to

(26) $$\sum_{i=1}^k q_{(i)} = A_1(\lambda_1) \quad \text{and} \quad \frac{\sum_{i=1}^k q_{(i)} Z_i}{\sum_{i=1}^k q_{(i)}} = \frac{A_2(\lambda_1)}{A_1(\lambda_1)}.$$

The second identity in (26) is

$$(27) \quad \frac{\sum_{i=1}^k \exp\{-\lambda_2 A_1(\lambda_1) Z_i \log(x_p/X_{n,n-k})/A_2^2(\lambda_1)\} Z_i}{\sum_{i=1}^k \exp\{-\lambda_2 A_1(\lambda_1) Z_i \log(x_p/X_{n,n-k})/A_2^2(\lambda_1)\}} = \frac{\log(x_p/X_{n,n-k})}{\log(A_1(\lambda_1)/p_n)}.$$

In order to demonstrate the existence of a solution to equation (11), first we will show that, with probability tending to one, there exists a continuous function $\lambda_2 = \lambda_2(\lambda_1)$ such that, for each $\lambda_1$, $(\lambda_1, \lambda_2) = (\lambda_1, \lambda_2(\lambda_1))$ is the solution to (27), and then we should prove that, for some $\lambda_1$, $(\lambda_1, \lambda_2) = (\lambda_1, \lambda_2(\lambda_1))$ is also the solution to the first identity of (26), and the solution satisfies both (12) and (13). To this end, set

$$f(\lambda) = \frac{\sum_{i=1}^k \exp\{-\lambda Z_i\} Z_i}{\sum_{i=1}^k \exp\{-\lambda Z_i\}}.$$

Then it is easy to see that $\lim_{\lambda \to -\infty} f(\lambda) = Z_1$, $\lim_{\lambda \to \infty} f(\lambda) = Z_k$ and $f(\lambda)$ is decreasing in $\lambda$ by checking that $\frac{d}{d\lambda} \log f(\lambda) < 0$. Therefore, there exists a unique continuous function $r(x)$ such that $f(r(x)) = x$ for any $x \in (Z_k, Z_1)$.

From now on we restrict $\lambda_1$ such that (12) holds, that is,

$$(28) \quad \frac{k}{n}\left(1 - \frac{\sqrt{\log(k/(np_n))}}{\sqrt{k}}\right) \leq A_1(\lambda_1) \leq \frac{k}{n}\left(1 + \frac{\sqrt{\log(k/(np_n))}}{\sqrt{k}}\right),$$

which implies

$$(29) \quad \left|\frac{\log(A_1(\lambda_1)/p_n)}{\log(k/(np_n))} - 1\right| \leq \frac{1}{\sqrt{k}\sqrt{\log(k/(np_n))}} \to 0$$

and

$$(30) \quad \begin{aligned} \frac{\log(x_p/X_{n,n-k})}{\log(k/(np_n))} &\left(1 - \frac{1}{\sqrt{k}\sqrt{\log(k/(np_n))}}\right) \\ &\leq \frac{\log(x_p/X_{n,n-k})}{\log(A_1(\lambda_1)/p_n)} \\ &\leq \frac{\log(x_p/X_{n,n-k})}{\log(k/(np_n))}\left(1 + \frac{1}{\sqrt{k}\sqrt{\log(k/(np_n))}}\right). \end{aligned}$$

Set $\mathcal{F}_n = \{Z_k < \frac{\log(x_p/X_{n,n-k})}{\log(k/(np_n))} < Z_1\}$. Then $P(\mathcal{F}_n) \to 1$ since

$$(31) \quad \frac{\log(x_p/X_{n,n-k})}{\log(k/(np_n))} \xrightarrow{p} \gamma^{-1}, \qquad Z_1 \xrightarrow{p} \infty, \qquad Z_k \xrightarrow{p} 0.$$

By definition, $f(r(\frac{\log(x_p/X_{n,n-k})}{\log(A_1(\lambda_1)/p_n)})) = \frac{\log(x_p/X_{n,n-k})}{\log(A_1(\lambda_1)/p_n)}$ on $\mathcal{F}_n$. This implies that

$$(32) \quad \lambda_2 = \lambda_2(\lambda_1) = r\left(\frac{\log(x_p/X_{n,n-k})}{\log(A_1(\lambda_1)/p_n)}\right) \frac{A_2^2(\lambda_1)}{A_1(\lambda_1) \log(x_p/X_{n,n-k})}$$

is the unique solution to equation (27), with probability tending to one.

Set
$$R_1 = \frac{\log(x_p/X_{n,n-k})}{\log(A_1(\lambda_1)/p_n)} - \frac{\log(x_p/X_{n,n-k})}{\log(k/(np_n))} \quad \text{and} \quad \hat{\lambda} = r\left(\frac{\log(x_p/X_{n,n-k})}{\log(A_1(\lambda_1)/p_n)}\right).$$

It follows from (30) that

(33) $$R_1 = O_p(k^{-1/2}(\log(k/(np_n)))^{-1/2})$$

holds uniformly for $\lambda_1$ under the restriction (12) or, equivalently, (28). Hereafter all terms $O_p(\cdot)$ and $o_p(\cdot)$ are assumed to hold uniformly for $\lambda_1$ if $\lambda_1$ is involved.

Using (6), (33) and

(34) $$\frac{\log(x_p/X_{n,n-k})}{\log(A_1(\lambda_1)/p_n)} - \hat{\gamma}_n^{-1} = \frac{\log(x_p/X_{n,n-k})}{\log(k/(np_n))} - \hat{\gamma}_n^{-1} + R_1$$
$$= -\frac{1}{\log(k/(np_n))} \log \frac{\hat{x}_p}{x_p} + R_1,$$

we have

(35) $$\frac{\log(x_p/X_{n,n-k})}{\log(A_1(\lambda_1)/p_n)} - \hat{\gamma}_n^{-1} = O_p(k^{-1/2}).$$

On the other hand, from Taylor's expansion, we have $f(\pm k^{-1/4}) - f(0) = \pm k^{-1/4} f'(0)(1 + o_p(1))$, where $f(0) = \hat{\gamma}_n^{-1}$, and

(36) $$f'(0) = -\left(\frac{1}{k}\sum_{i=1}^{k} Z_i^2 - \hat{\gamma}_n^{-2}\right) = -\gamma^{-2}(1 + O_p(k^{-1/2})).$$

For the proof for the last step of (36), see, for example, [4] or [9]. Hence,
$$P\left(f(k^{-1/4}) < \frac{\log(x_p/X_{n,n-k})}{\log(A_1(\lambda_1)/p_n)}\right) \to 1$$

and
$$P\left(f(-k^{-1/4}) > \frac{\log(x_p/X_{n,n-k})}{\log(A_1(\lambda_1)/p_n)}\right) \to 1.$$

Therefore, $\hat{\lambda} = r(\frac{\log(x_p/X_{n,n-k})}{\log(A_1(\lambda_1)/p_n)})$ satisfies $P(\hat{\lambda} \in (-k^{-1/4}, k^{-1/4})) \to 1$. Thus, from Taylor's expansion, we obtain
$$f(\hat{\lambda}) - \hat{\gamma}_n^{-1} = \hat{\lambda} f'(0)(1 + O_p(k^{-1/4})) = -\gamma^{-2}\hat{\lambda}(1 + O_p(k^{-1/4})),$$

which, coupled with (34), yields

(37) $$\hat{\lambda} = -\gamma^2 \left(\frac{\log(x_p/X_{n,n-k})}{\log(A_1(\lambda_1)/p_n)} - \hat{\gamma}_n^{-1}\right)(1 + O_p(k^{-1/4}))$$
$$= \frac{\gamma^2}{\log(k/(np_n))} \log \frac{\hat{x}_p}{x_p}(1 + O_p(k^{-1/4})).$$

Then, by using (35) and Taylor's expansion, we have

$$\text{(38)} \qquad \sum_{i=1}^{k} \exp\{-\hat{\lambda} Z_i\} = k(1 + O_p(k^{-1/2})).$$

Note that

$$\lambda_2 = \lambda_2(\lambda_1) = \frac{\hat{\lambda} A_2^2(\lambda_1)}{A_1(\lambda_1) \log(x_p/X_{n,n-k})},$$

where $\hat{\lambda}$ is a function of $\lambda_1$ as well. Plug $\lambda_2 = \lambda_2(\lambda_1)$ into $q_{(i)}$ and set $h(\lambda_1) = \sum_{i=1}^{k} q_{(i)} - A_1(\lambda_1)$. Then $h(\lambda_1)$ has the expression

$$\frac{1}{n} \exp\{-1 - \lambda_1\} \exp\left\{ \frac{\hat{\lambda} A_2^2(\lambda_1)}{A_1(\lambda_1) \log(x_p/X_{n,n-k})} \right.$$

$$\left. \times \left( \frac{\log(x_p/X_{n,n-k})}{A_2(\lambda_1)} - \frac{1}{A_1(\lambda_1)} \right) \right\} \sum_{i=1}^{k} \exp\{-\hat{\lambda} Z_i\} - A_1(\lambda_1).$$

Put

$$\lambda_1' = -\log\left(1 + \frac{\sqrt{k}\sqrt{\log(k/(np_n))}}{n-k}\right) - 1$$

and

$$\lambda_1'' = -\log\left(1 - \frac{\sqrt{k}\sqrt{\log(k/(np_n))}}{n-k}\right) - 1.$$

It is easy to check that

$$e^{-1-\lambda_1'} = 1 + \frac{\sqrt{k}\sqrt{\log(k/(np_n))}}{n-k},$$

$$\text{(39)} \qquad A_1(\lambda_1') = \frac{k}{n}\left(1 - \frac{\sqrt{\log(k/(np_n))}}{\sqrt{k}}\right),$$

$$A_2(\lambda_1') = \gamma^{-1}\frac{k}{n}\left(1 - \frac{\sqrt{\log(k/(np_n))}}{\sqrt{k}} + o_p\left(\frac{\sqrt{\log(k/(np_n))}}{\sqrt{k}}\right)\right).$$

The first two identities are obvious. The third one follows from the second one, equation (35) and a well-known result for the Hill estimator, that is, $\sqrt{k}(\hat{\gamma}_n^{-1} - \gamma^{-1}) \xrightarrow{d} N(0, \gamma^{-2})$.

Now it follows from (39) and (38) that

$$h(\lambda_1') = \frac{k}{n} \frac{\sqrt{\log(k/(np_n))}}{\sqrt{k}}(1 + o_p(1)).$$

Similarly,
$$h(\lambda_1'') = -\frac{k}{n}\frac{\sqrt{\log(k/(np_n))}}{\sqrt{k}}(1+o_p(1)).$$

Hence, with probability tending to one, there exists a $\lambda_1$ satisfying (12) and the first equation in (26) with $\lambda_2 = \lambda_2(\lambda_1)$ defined in (32); that is, we have shown the existence of the solution to (11) such that (12) and (13) hold.

We still need to estimate $\hat{\lambda}_1$ from the equation $h(\hat{\lambda}_1) = 0$, which is equivalent to

(40)
$$\frac{1}{n}\exp\left\{\hat{\lambda}\left(\frac{\log(x_p/X_{n,n-k})}{\log(A_1(\hat{\lambda}_1)/p_n)} - \frac{\log(x_p/X_{n,n-k})}{(\log(A_1(\hat{\lambda}_1)/p_n))^2}\right)\right\}\sum_{i=1}^{k}\exp\{-\hat{\lambda}Z_i\}$$
$$= \exp\{1+\hat{\lambda}_1\} - \frac{n-k}{n}.$$

It follows from (35) that
$$\frac{\log(x_p/X_{n,n-k})}{\log(A_1(\hat{\lambda}_1)/p_n)} - \frac{\log(x_p/X_{n,n-k})}{(\log(A_1(\hat{\lambda}_1)/p_n))^2} = \hat{\gamma}_n^{-1} + o_p(1).$$

Hence, applying Taylor's expansion to both sides of (40) yields

(41) $$1+\hat{\lambda}_1 = O_p\left(\frac{\sqrt{k}}{n}\right).$$

It is easy to show that $\max_{1\le i\le n}|nq_i - 1| = o_p(1)$. Thus,
$$L(x_p) = 2\sum_{i=1}^{n} nq_i\{nq_i - 1 - \tfrac{1}{2}(nq_i-1)^2(1+o_p(1))\}$$
$$= \sum_{i=1}^{n}(nq_i - 1)^2(1+o_p(1))$$
$$= \sum_{i=1}^{n}(\log(nq_i))^2(1+o_p(1))$$
$$= (1+o_p(1))\left((n-k)(1+\hat{\lambda}_1)^2 + \sum_{i=1}^{k}(\log(nq_{(i)}))^2\right).$$

It follows from (29) and (35) that
$$\log(nq_{(i)}) = -(1+\hat{\lambda}_1) + \hat{\lambda}\left(\frac{\log(x_p/X_{n,n-k})}{\log(A_1(\hat{\lambda}_1)/p_n)} - \frac{\log(x_p/X_{n,n-k})}{(\log(A_1(\hat{\lambda}_1)/p_n))^2} - Z_i\right)$$
$$= -(1+\hat{\lambda}_1)$$
$$\quad - \hat{\lambda}\left(Z_i - \hat{\gamma}_n^{-1} + O_p\left(\frac{\sqrt{\log(k/(np_n))}}{\sqrt{k}} + \frac{1}{\log(k/(np_n))}\right)\right),$$

uniformly for $i = 1, \ldots, k$. As in (36), we have

$$\sum_{i=1}^k (Z_i - \hat{\gamma}_n^{-1})^2 = k\left(\frac{1}{k}\sum_{i=1}^k Z_i^2 - \hat{\gamma}_n^{-2}\right) = k\gamma^{-2}(1 + O_p(k^{-1/2})).$$

Hence, by (6), (37) and (41),

$$\begin{aligned}
L(x_p) &= (1 + o_p(1))\left(n(1 + \hat{\lambda}_1)^2 + \hat{\lambda}^2 \sum_{i=1}^k (Z_i - \hat{\gamma}_n^{-1})^2\right) \\
&= \gamma^2 k \left\{\frac{\log(\hat{x}_p/x_p)}{\log(k/(np_n))}\right\}^2 + o_p(1) \\
&\xrightarrow{d} \chi_1^2.
\end{aligned}$$ $\square$

**Acknowledgments.** The authors would like to thank Professor Richard Green for his careful reading of the paper and his comments. The comments from an Editor, an Associate Editor and three reviewers were very helpful.

School of Mathematics  
Georgia Institute of Technology  
Atlanta, Georgia 30332-0160  
USA  
E-mail: [peng@math.gatech.edu](peng@math.gatech.edu)

Department of Mathematics  
and Statistics  
University of Minnesota Duluth  
SCC 140, 1117 University Drive  
Duluth, Minnesota 55812  
USA  
E-mail: [yqi@d.umn.edu](yqi@d.umn.edu)